\documentclass[11pt]{article}

\newtheorem{thm}{Theorem}[section]

\newtheorem{prop}[thm]{Proposition}
\newtheorem{lem}[thm]{Lemma}
\newtheorem{cor}[thm]{Corollary}

\newtheorem{defn}[thm]{Definition}

\newtheorem{example}[thm]{Example}
\newtheorem{remark}[thm]{Remark}
\newtheorem{conj}[thm]{Conjecture}

\newcommand{\PPP}{{\bf P}}
\newcommand{\PP}[1]{{{\bf P}^{#1}}}
\newcommand{\pr}[1]{{{\bf P}^{#1}}}
\newcommand{\skipit}[1]{{}}
\newcommand{\prfend}{\hbox to7pt{\hfil}
\par\vskip-\baselineskip\hbox to\hsize
{\hfil\vbox {\hrule width6pt height6pt}}\vskip\baselineskip}

\newcommand{\myarrow}[2]{\hbox to #1pt{\hfil$\to$\hfil}{\hskip-#1pt{\raise
10pt\hbox to#1pt{\hfil$\scriptscriptstyle #2$\hfil}}}}

\def\a{\bigskip \par \noindent}
\def\b{\par \noindent}

\begin{document}
\title{The role of the cotangent bundle in resolving ideals of fat points in the
plane.}

\author{Alessandro Gimigliano\\
\and
Brian Harbourne\\
\and
Monica Id\`a\\}

\maketitle

\thanks{Acknowledgments: We thank GNSAGA, MUR and
the University of Bologna, which supported
visits to Bologna by the second author, who also
thanks the NSA and NSF for supporting his research.}

\begin{abstract}
We study the connection between the generation of a fat point scheme supported
at general points
in $\PP2$ and the behaviour of the cotangent bundle with respect to some
rational curves particularly relevant for the scheme. We put forward two
conjectures, giving examples and partial results in support of them.

\end{abstract}

\section{Introduction}\label{intro}

In this paper we are concerned with minimal free graded resolutions
of fat point ideals in $\PP2$. Given general points
$P_1,\ldots,P_n\in \pr2$ (which, unless we say something explicit
to the contrary, will always be assumed
to be general), and nonnegative integers $m_1,\ldots,m_n$,
let $I(Z)$ denote the ideal $I(P_1)^{m_1}\cap \cdots I(P_n)^{m_n}$
of $R=K[\pr2]=K[x_0,x_1,x_2]$ (where $K$ is any algebraically closed field
and where $I(P_i)$ is the ideal generated by all forms that vanish at $P_i$).
We refer to $I(Z)$ as a {\it fat point\/} ideal, and if $Z$ is
the subscheme defined by $I(Z)$, we use
  $m_1P_1+\cdots+m_nP_n$ or $Z(m_1,\ldots,m_n)$ to denote the scheme $Z$, and
${\cal I}_{Z}$ for its sheaf of ideals, so that, in particular,
$I(Z)_k=H^0(\PP2,{\cal I}_{Z}(k))$.

In order to understand better the geometry of $Z$ as a subscheme of $\PP2$,  the
first thing that comes to mind is to see how many curves of given degree $k$
contain $Z$, that is, have singularities of multiplicity at least
$m_1,\ldots,m_n$ at the given points $P_1,\ldots,P_n$; in other words, we want
to determine the dimension, as a $K$-vector space, of the homogeneous component
$I(Z)_k$
of $I(Z)$.
\par The Hilbert function $\, h_Z$ of $\, I(Z)$, $\, h_Z(k) :=
\hbox{dim}_K(I(Z)_k)$, is not known in general, even if it has been determined
for many choices of $Z$. For example, it is known for
all $Z$ with $n\le 9$ (\cite{refNtwo}, or see \cite{refanti}), for any $n$ if
$m_1=\cdots=m_n\le 20$ (\cite{refHi1}, \cite{refCCMO}), and for any $n$ if
$m_i\le 7$ (\cite{refmig}, \cite{refyang}).
(Some but not all of these and later citations assume $K$
is the complex numbers.)
It is also known for many additional cases.
Let us say that the sequence of multiplicities $m_i$
(and by extension $Z$) is {\it uniform}
if $m_1=\cdots=m_n\ge0$ and if $n\ge 9$.
Then \cite{refEvain} determines $h_Z$ for all $k$,
as long as $Z$ is uniform, if $n$ is a square, extending
results of \cite{refHHF}.
The paper \cite{refHRa} determines $h_Z$ in
many other uniform cases.
\par All of these results are consistent
with a well known conjecture by means of which
one can explicitly write down the function
$h_Z$ given the multiplicities $m_i$.
Various equivalent versions of this conjecture have been given
(see \cite{refSe}, \cite{refVanc}, \cite{refGi}, \cite{refHi},
\cite{refSiena}). We will refer to them collectively as the
SHGH Conjecture.
\par Let us say that a fat point subscheme $Z$ is {\it quasi-uniform\/}
if $n\ge 9$ and $m_1=\cdots=m_{9}\ge m_{10}\ge \cdots \ge m_n\ge0$.
Thus uniform implies quasi-uniform. As shown in \cite{refHHF}, assuming the
SHGH Conjecture, then
$h_Z(k)= \hbox{max}(0, {k+2\choose 2}-\sum_i{m_i+1\choose2})$
holds for all $k$ for a quasi-uniform $Z$. Since there are
${k+2\choose 2}$ forms of degree $k$ and since the requirement
for a form to vanish to order $m_i$ at a point $P_i$ imposes
${m_i+1\choose 2}$ conditions, the SHGH Conjecture in this situation
just says that the conditions imposed by the points are independent
as long as $h_Z(k)>0$.

\par \medskip To go deeper into the geometry of a fat point scheme, the next
step consists in understanding the relations among the curves containing $Z$,
that is, determining
the minimal free graded resolution $0\to M_1\to M_0\to I(Z)\to 0$
of $I(Z)$. Here $M_0$ and $M_1$ are free $R$-modules of the form $M_0=\oplus_k
R^{t_k}[-k]$ and
$M_1=\oplus_k R^{s_k}[-k]$. If $h_Z$ is known and if the graded Betti numbers
$t_k$ are known, then the values of $s_k$ are easy to determine from the exact
sequence above.
\par We are hence interested in the graded Betti numbers $t_k$.
It is not hard to see that $t_k$ is
the dimension of the cokernel of the map
$\mu_{k-1}(Z): I(Z)_{k-1}\otimes R_1\to I(Z)_k$, where $R_1$ denotes
the $K$-vector
space spanned in $R$ by linear forms and $\mu_{k-1}$
is the map induced by multiplication of elements of $I(Z)_{k-1}$
by linear forms.
This paper is a reflection about the geometric obstacles to the rank
maximality of the maps $\mu_k$.  Let us denote by $\Omega$ the cotangent
bundle of $\PP2$, and by $p:X\to \PP 2$ the blow up at the points $P_i$. We
first translate the problem of determining the rank of the maps $\mu_{k}(Z)$
into two equivalent postulation problems for $Z$, one in $\PP 2$ and the other
in $X$: determine, for each $k$, the rank of the restriction map
\par (a) ${\rho_k}={\rho_k}(Z): H^0(\Omega(k+1)) \to H^0(\Omega(k+1)|_Z)$; or
\par (b) ${\eta_k}={\eta_k}(Z): H^0(p^*\Omega(k+1)) \to
H^0(p^*\Omega(k+1)|_{p^{-1}Z})$.
\b We show that the point of view (a) gives some information about the failure
of this rank maximality due to superfluous conditions imposed by $Z$ to the
restriction of $\Omega$ to some curves; but in fact this is not enough, and the
right point of view is (b), since it is then possible to take into account the
splitting of $p^*\Omega$ on the normalization of the appropriate rational
curves, and this allows to count properly the superfluous conditions imposed by
$Z$ to the restriction of $\Omega$ to each curve.

Hence, by studying several examples and proving certain results (e.g.
\ref{betarhocor}), we arrive at two conjectures about the failure of
the rank maximality of $\mu _k$, one when $\mu _k$ is expected
to be surjective and the other when injectivity is expected.
 The idea is the same in the two cases but the expected surjective case is much
easier to formulate, and this is why we keep them distinct; in both cases the
obstruction to rank maximality is described by the presence of particular
rational curves whose intersection with $Z$ is ``too high".

\par \medskip  Notice that a similar line of thought leads to the SHGH
Conjecture: in fact, determining $h^0(\PP2, {\cal I}_Z(k))$ amounts to
computing the rank of the restriction map
$r_k: H^0(\PP2, {\cal O}_{\PP2}(k)) \to H^0(Z,{\cal O}_Z)$.
The SHGH Conjecture says that failure of $r_k$
to have maximal rank is completely accounted for by
the occurrence of curves $C\subset \PP2$ whose strict transform
$\tilde C\subset X$
is an exceptional divisor (i.e., a smooth rational curve of
self-intersection $-1$),
such that the scheme-theoretic intersection $C\cap Z$ is too big with
respect to
${\cal O}(k)|_C$ (or, expressing things on the blow up, such that the inverse
image $\tilde Z$ of $Z$ meets $\tilde C$ in too many points with respect to
$kL|_{\tilde C}$, which here just means that $\tilde C\cdot F<-1$,
where $F=kL-m_1E_1-\cdots -m_nE_n$, $L$ is
the pullback to $X$ of a general line in $\PP2$ and $E_i$ is the exceptional
locus
obtained by blowing up the point $P_i$).

\par Unfortunately, things are quite complicated when studying the postulation
with respect to a rank 2 vector bundle; for example, as said above, we have to
take into consideration the splitting of $p^*\Omega$ on the normalization of a
rational plane curve, which is not known in general (see \ref {splitlem}), and
is actually an interesting problem {\em per se}. In our examples we have made
use, when necessary, of a Macaulay 2 script which allows us to compute
splitting types (see Section A2.3 of \cite{refGHIpv}).

\par The use of the cotangent bundle in problems concerning the generation of
homogeneous ideals of subschemes of
a projective space was introduced by A.Hirschowitz, and used for the first time
for curves in $\PP3$ (see \cite{refId3}).

\par Our conjectures assume that the fat point scheme $Z$ postulates well in the
degree $k$ we are considering, i.e. that $h^1({\cal I}_{Z}(k))=0$; but notice
that, assuming the SHGH Conjecture, we can always reduce ourselves to
considering fat points $Z$ with good postulation, and for these we need to
study only the map $\mu_{\alpha}$, where $\alpha$ is the initial degree of
$I(Z)$ (see \ref {CastelMumRegRem}).

\par \medskip Here is what is currently known about resolution of fat point
ideals in $\PP2$. For uniform (\cite{refIGP})
or quasi-uniform (\cite{refHHF}) $Z$, it is conjectured that the maps $\mu_k$
have maximal rank for all $k$. We refer to these as the Uniform Resolution
and Quasi-Uniform Resolution Conjectures.
The Uniform Resolution Conjecture has been proved for $m=1$ (\cite{refGM}),
$m=2\,$ (\cite{refId1}) and $m=3\,$ (\cite{refGiId}); more generally, if
$m_i\le3$ for all $i$, and the length of $Z$ is sufficiently high,
\cite{refBI} determines the graded Betti numbers in all degrees. Verifications
of the
Quasi-Uniform Resolution Conjecture in some
cases were given in \cite{refHHF},
under the assumption of the SHGH Conjecture.
Some outright verifications were given by \cite{refHRa}.
By applying the results of \cite{refEvain} to results of \cite{refHHF},
it also follows that the Uniform Resolution Conjecture
holds for all $m$ not too small, as long as $n$ is
an even square. Finally, the Betti numbers are known
for all $Z$ with $n\le8$ (\cite{refCa}, \cite{refF3}, \cite{refCJM},
\cite{refFHH}); the $n\leq8$ results
show that any general resolution conjecture
will have to be more subtle than the SHGH Conjecture.
\par In Section \ref{twoconj} we prove that our
Conjectures \ref{conj1} and \ref{conj2} together with the SHGH Conjecture
imply the Uniform and Quasi-uniform Resolution Conjectures (see Proposition
\ref{conjqunif}).

\section{Preliminaries}\label{prelims}

We now establish some terminology and notations and
recall some basic concepts.
\par \bigskip By curve we will mean a 1-dimensional scheme without embedded
components.

The surface obtained from $\PP2$ by blowing
up general points $P_i$ is always denoted by $X$, $p: X\to \pr2$ is
the morphism
given by blowing up the points, $E_i$ is the exceptional curve
obtained by blowing up the point $P_i$ and $L$ is
the divisorial inverse image under $p$ of a line in $\PP2$.
We will also use $L$ and $E_i$ to denote the linear equivalence
class of the given divisor, in which case
the divisor class group $\hbox{Cl}(X)$ is the free abelian group
on the basis $L, E_1,\ldots, E_n$.
The intersection form on $X$ is such that the basis elements
are orthogonal with $-L^2=E_i^2=-1$ for all $i$.
\par Given a divisor $F$ on $X$, we will use $F$ to denote its divisor class
and sometimes even the sheaf ${\cal O}_X(F)$, and
we will for convenience write $H^0(F)$ for $H^0(X, {\cal O}_X(F))$. For each
$F$,
there is a natural multiplication map
$\mu_F: H^0(F)\otimes H^0(L) \to H^0(F+L)$.

If $Z=m_1P_1+\cdots+m_nP_n$ is a fat point scheme, it is clear that, under the
correspondence of $H^0(\PP2,{\cal I}_{Z}(k))$ with
$H^0(X, kL -\sum m_iE_i)$, the map
$$\mu_k(Z): H^0({\cal I}_{Z}(k))\otimes H^0({\cal O}_{\PP 2}(1))\to H^0({\cal
I}_{Z}(k+1))$$
is just the map $\mu_{tL-\sum m_iE_i}$.

Given a curve $C\subset\PP 2$, we denote the multiplicity of $C$ at $P_i$
by $m(C)_{P_i}=r_i$, and $\tilde C= dL-\sum r_iE_i$ will denote its
strict transform.
Note that $d$ is just the degree of $C$.
If $C\subset \PP 2$ is an integral curve such that
$\tilde C\subset X$ is smooth and rational, we write ${\cal O}_{\tilde C}(k)$
instead of ${\cal O}_{\PP 1}(k)$.
We recall that $\tilde C$ is an exceptional divisor (of the first kind)
in $X$ if $\tilde C= dL-\sum r_iE_i$ is smooth and rational with $-1=\tilde
C^2=d^2-\sum r_i^2$, which
by the adjunction formula implies $-1=K_X\cdot \tilde C=-3d+\sum r_i$, since
$K_X=-3L+E_1+\cdots+E_n$.

\medskip Let $Y$ be a smooth projective variety, $D$ a divisor and $A$ a
subscheme  of $Y$; the residual scheme $A'=\hbox{res}_DA$ is the subscheme of
$Y$ whose sheaf of ideals ${\cal I}_{\hbox{res}_DA}$ is given by the exact
sequence: $0\to {\cal I}_{\hbox{res}_DA}(-D) \to {\cal I}_{A}\to {\cal
I}_{A\cap D,D}\to 0$, where ${\cal I}_{A\cap D,D}$ is the sheaf of ideals on
$D$ defining the scheme-theoretic intersection of $A$ and $D$ as a subscheme of
$D$.

\par If $Z=m_1P_1+\cdots+m_nP_n$ in $\PP2$  is a fat point scheme, and $C$ is a
plane curve
whose proper transform is $\tilde C=dL-\sum r_iE_i$, the residual sequence
tensored by ${\cal O}_{\PP2}(k)$ becomes: $0\to {\cal I}_{Z'}(k-d) \to
{\cal I}_{Z}(k)\to {\cal I}_{Z\cap C,C}(k)\to 0$, where  $Z'=\hbox{res}_CZ$
has homogeneous ideal $(I(Z):I(C))$.
\b Now if we set $F_k(Z)=kL-m_1E_1-\cdots-m_nE_n$, we have $F_k(Z)-{\tilde
C}=(k-d)L-\sum(m_i-r_i)E_i$ and its cohomology is the cohomology of a fat point
scheme provided that $m_i-r_i\geq 0$ for all $i$; more precisely,
$F_k(Z)-{\tilde C}=F_{k-d}(Z')$ if $r_i\leq m_i$. Thus divisors corresponding
to residuals are easy to compute.

\par \bigskip Setting $\Omega=\Omega_{\PP 2}$, recall the Euler sequence on $\PP
2$:
$$0 \to \Omega (1) \to {\cal O}_{\PP 2} \otimes H^0({\cal O}_{\PP 2}(1)) \to
{\cal O}_{\PP 2}(1) \to 0.$$

Now let $C\subset \PP 2$ be a degree $d$ integral curve and assume
$\tilde C\subset X$ smooth and rational. Since the Euler sequence is a sequence
of vector bundles, its pullback to $X$
 restricted to $\tilde C$ is still exact,
and gives
$$0 \to p^*\Omega (1)|_{\tilde C} \to
{\cal O}_{\tilde C}\otimes H^0({\cal O}_{X}(L)) \to {\cal O}_{\tilde
C}(d) \to 0.
\eqno(*)$$
In the following we set
$$p^*\Omega(1)|_{\tilde C}\cong{\cal O}_{\tilde C}(-a_C)\oplus{\cal
O}_{\tilde C}(-b_C),$$
where we always assume $a_C\leq b_C$; looking at  the Chern classes in $(*)$
gives $a_C+ b_C=d$. We will say that the splitting {\em type}
of $C$ or $\tilde C$ is $(a_C,b_C)$ and the splitting {\em gap} is $b_C-a_C$.

In some cases we can immediately determine the splitting type.
Suppose that $m$ is the maximum
value of $m(C)_{P_i}$.
See \cite{refAs} or \cite{refF1}, \cite{refF2} for the proof of
the following lemma:

\begin{lem}\label{splitlem}
We have $\hbox{min}(m,d-m)\le a_C\le d-m$, and $d=a_C+b_C$.
\end{lem}

Note that the splitting type is completely determined if
$d-m\le m+1$, and it is $(\hbox{min}(m,d-m),\hbox{max}(m,d-m))$. When $d-m>m+1$
it is not known in general
what the splitting type is, but it
can be computed fairly efficiently; see Section A 2.3 in \cite{refGHIpv}.

\par \bigskip If $f:A \to B$ is a linear map between vector spaces, we say that
$f$
is {\em exp-onto} (i.e., expected to be onto), resp. {\em exp-inj} (i.e.,
expected to be injective), if $\hbox{dim }A \geq \hbox{dim }B$,
resp. $\hbox{dim }A \leq \hbox{dim }B$.
The {\em expected dimension} for the cokernel of $f$ is defined to be
$\hbox{exp-dim cok}(f):= \hbox{max}(0, \hbox{dim }B - \hbox{dim }A)$.
So, for example,
$$\hbox{exp-dim cok}(\mu_k(Z))=\hbox{max}(0, h^0({\cal
I}_Z(k+1))-3h^0({\cal I}_Z(k))).$$

We say that a fat point scheme
$Z$ has {\em good postulation in degree} $k$, if the map $r_k$
is of maximal rank, i.e. if $h^0({\cal I}_Z(k))h^1({\cal I}_Z(k))=0$.
We say that $Z$ has {\em good postulation} if the maps $r_k$ have maximal rank
for all $k$, and we say that $Z$ is {\em minimally generated} if
the maps $\mu_k$ all have maximal rank (i.e., $Z$ is
minimally generated if $\mu_k$ is onto when it is exp-onto
and injective when it is exp-inj).

\medskip A few additional notions will be useful.
Given a 0-dimensional scheme $Y$, we denote by $l(Y)$ the length of $Y$;
hence $l(m_1P_1+\cdots+m_nP_n)=\sum_i {m_i+1 \choose 2}$.
\par We define $\alpha=\alpha(Z)$ to be the least
$k$ such that $h^0({\cal I}_Z(k))$ is positive, and
we define $\tau=\tau(Z)$ to be the least $k$ such that
$h^1({\cal I}_Z(k))=0$.
\smallskip \b Recall that if $h^1({\cal I}_Z(k))=0$ then $h^1({\cal I}_Z(t))=0$
for $t\geq k$, and $\mu_t(Z)$ is surjective for $t\geq k+1$, by the
Castelnuovo-Mumford lemma
\cite{refMu2}.

\begin{remark}\label{numinvlem}
\rm Let $Z$ be a fat points subscheme of $\PP2$ (supported at general points).
Then $\alpha-1\leq\tau$. If $Z$ has good postulation, then
$\alpha-1\leq\tau\leq\alpha$.

\par In fact, $\alpha-1\leq\tau$
follows by taking cohomology of
$0\to {\cal I}_Z(k) \to {\cal O}_{\PP2}(k) \to {\cal O}_Z\to 0$.
Good postulation gives $h^0({\cal I}_Z(k))h^1({\cal I}_Z(k))=0$,
which implies $\tau\leq\alpha$.
\end{remark}

\begin{remark}\label{CastelMumRegRem} \rm
Since $\mu_k(Z)$ (being the 0-map) is trivially injective for all $k<\alpha$
and it is surjective for $k\geq \tau+1$, we need only consider
$\mu_k$ in degrees $k$ (if any) with $\alpha\leq k\leq \tau$.

If $Z$ has good postulation, then either $\tau = \alpha-1$, and the Betti
numbers for $I(Z)$ are completely determined, or $\tau =\alpha$, in which case
we need only consider $\mu _{\alpha}$; if $\mu _{\alpha}$ is exp-onto, then
$Z$ is minimally generated if and only if $\mu _{\alpha}$ is surjective,
while if $\mu _{\alpha}$ is exp-inj,
$Z$ is minimally generated if and only if $\mu _{\alpha}$ is injective.
\par Now drop the good postulation assumption, and take any $Z$; if
$k\ge\alpha$, assuming the SHGH Conjecture
it is always possible (and easy to do explicitly,
by factoring out the fixed part of $H^0({\cal I}_Z(k))$;
see \cite{refGH}) to replace $k$ and $Z$ by a $k'$ and $Z'$ (supported at the
same points)
such that the kernels of $\mu_k(Z)$ and $\mu_{k'}(Z')$ have the same
dimension, but such that $Z'$ has good postulation in degree $k'$.
Thus (assuming the SHGH Conjecture) we can reduce to considering only
fat points $Z$ with good postulation and with $\alpha = \tau$, and
for these we need to study only the map $\mu _{\alpha}$.
\end{remark}

\a The forthcoming Remark \ref{OmegaDimRem} and Lemma \ref{lengthlem} will be
useful in the next section:
\begin{remark}\label{OmegaDimRem} \rm
Let $C$ be a curve of degree $d$ in $\PP 2$; then the exact sequence
$0 \to {\cal O}_{\PP 2}(t-d) \to {\cal O}_{\PP 2}(t)
\to {\cal O}_{\PP 2}(t)\vert_{C}  \to 0$
gives
\a $h^0({\cal O}_{\PP 2}(t)\vert_{C}) = {t+2 \choose 2}$ for
$0 \leq t \leq d-1$,  $\; \; h^0({\cal O}_{\PP 2}(t)\vert_{C}) =
{1 \over 2}(2td+3d-d^2)$ for $t \geq d$.
\a The same exact sequence twisted
by $\Omega$ and the cohomology of the cotangent bundle (see for example
\cite{refOSS}):
$$h^0(\PP2, \Omega(k))=h^2(\PP2, \Omega(-k))=\left\{ \begin{array}{ll}
                           k^2-1 & \hbox{if $k\geq 1$} \\
                           0     & \hbox{if $k\leq 0$}
                                   \end{array}
                           \right. ,
\quad h^1(\PP2, \Omega(k))=\left\{ \begin{array}{ll}
                           0  & \hbox{if $k\neq 0$} \\
                           1  & \hbox{if $k=0$}
                                   \end{array}
                           \right. $$
together give:

$$h^0({\Omega}(t)\vert_C)=\left\{ \begin{array}{ll}
                           (t-1)(t+1) & \hbox{if $1 \leq t \leq d-2$} \\
                           d(2t-d)     & \hbox{if $t \geq d-1$}
                                   \end{array}
                           \right. ,
\quad h^1({\Omega}(t)\vert_C) = 0 \quad \hbox{for $t\geq \hbox{max}(1, d-1).$}
$$

\end{remark}

\begin{lem}\label{lengthlem}
Let $C$ be a plane curve having a singularity
of multiplicity $r$ at a point $P$, and let $Z$ be the $m$-fat point
supported at $P$; then $l(Z\cap C)={m+1\choose 2}-{m-r+1\choose 2}$.
\end{lem}

\noindent {\bf Proof.} Let $x, y$ be local coordinates at $P$,
$\,A_m:=K[x,y]/(x,y)^m$ the coordinate ring of $Z\,$, $f=0$ a local equation
for $C$, where $f$ has initial degree $r$, and $({\bar f}):=(f)A_m\,$; then,
$l(Z \cap C)$ is the dimension of the $K$-vector space $A_m {\big /}({\bar
f})$. If $r\geq m$,  ${\bar f}=0$,  so  ${\rm dim}_{K}A_m {\big
/}({\bar f})={m+1 \choose 2}$ (which was already obvious since $Z \subset C$).
If $r < m$,
it is easy to prove that the vector space $({\bar f})$ has dimension
${m+1-r\choose 2}$ using an appropriate induction.

\prfend

\section{Various equivalent postulation problems}\label{setup}

In this and in the following sections $k$ will always denote a positive integer,
and $Z$, as usual, a fat point subscheme
supported at general points of $\PP2$.
\par In this section we are going to translate the problem of determining the
rank of the maps $\mu_{k}(Z): I(Z)_{k}\otimes R_1\to I(Z)_{k+1}$ into three
different, but closely related, postulation problems. By {\em postulation
problem} we mean the computation of the rank of a restriction map $H^0(F) \to
H^0(F\vert_{Y})$ with $F$ a vector bundle and $Y$ a subscheme of a given
scheme.
One of these approaches, i.e. the translation into a postulation problem in the
3-fold $\PPP(\Omega)$ with respect to a rank 1 bundle, is here because we find
it intrinsically interesting, altough we'll use it only to understand the
geometry of certain examples. The other two approaches will lead to conjectures
\ref{conj1} and \ref{conj2}.

\a We now define the three restriction maps in which we are interested:
$\rho_k=\rho_k(Z)$, $\psi_k=\psi_k(Z)$, $\eta_k=\eta_k(Z)$.
\par  The multiplication map $\mu_k=\mu_k(Z)$ comes from
considering the Euler sequence twisted by ${\cal I}_Z(k)$ and taking cohomology:

\vskip6pt
\b $ (1*) \quad 0 \to H^0(\Omega (k+1)\otimes {\cal I}_Z)\to
H^0({\cal I}_Z(k))\otimes H^0({\cal O}_{\PP 2}(1))
\buildrel {\mu_k} \over \to
H^0({\cal I}_Z(k+1))\to H^1(\Omega (k+1)\otimes {\cal I}_Z)
\to H^1({\cal I}_Z(k))\otimes H^0({\cal O}_{\PP 2}(1))\to \dots $\vskip6pt

In the forthcoming Lemma \ref{postlem}
we compare this to
the cohomology sequence obtained by
restricting $\Omega$ to $Z$:

$$(2*)\quad  0 \to H^0(\Omega (k+1)\otimes {\cal I}_Z)\to
H^0(\Omega (k+1)) \buildrel {\rho_k} \over \to H^0(\Omega
(k+1)\vert_{Z})\to H^1(\Omega (k+1)\otimes {\cal I}_Z) \to
H^1(\Omega (k+1))=0  $$

Now consider the projective bundle
$\; \pi:\PPP(\Omega)\to \PP 2$ with the invertible sheaf
$${\cal E}_t={\cal O}_{\PPP(\Omega)}(1)\otimes \pi^* {\cal O}_{\PP 2}(t).$$
We set $$T= \pi^{-1}(Z)\subset \PPP(\Omega).$$

By \cite{refHt} Ex. III.8.1, III.8.3 and III.8.4,
$R^i\pi_*{\cal O}_{\PPP(\Omega)}(1)=0$
for $i > 0$, hence
$R^i\pi_*{\cal E}_t \cong R^i\pi_*{\cal O}_{\PPP (\Omega)}(1)\otimes
{\cal O}_{\PP 2}(t)=0$
for $i > 0$, so that $H^i(\Omega (t)) \cong H^i({\cal E}_t)$ for
all $i \geq 0$; in particular, $H^1({\cal E}_{k+1})=0$ for any $k \geq 0$.
Taking $\psi_k$ to be the canonical restriction map, we have the exact sequence:

$$ (3*)\quad  0 \to H^0({\cal E}_{k+1}\otimes {\cal I}_{T})\to
H^0({\cal E}_{k+1}) \buildrel {\psi _k} \over \to H^0({\cal E}_{k+1\vert T})\to
H^1({\cal E}_{k+1}\otimes {\cal I}_{T}) \to 0.$$

We will also work in the blow up $p:X \to \PP 2$.
Set $$\tilde Z= \sum_{i\geq 1} m_iE_i \subset X$$
 and consider the exact sequence:

$$(4*) \qquad 0 \to H^0(p^*\Omega (k+1)\otimes
{\cal I}_{\tilde Z})\to H^0(p^*\Omega (k+1)) \buildrel {\eta _k}
\over \to H^0(p^*\Omega (k+1)\vert_{{\tilde Z}})\to
H^1(p^*\Omega (k+1)\otimes {\cal I}_{\tilde Z}) \to 0$$

\noindent where $H^1(p^*\Omega (k+1))=0$ for the following reason:
$R^ip_*{\cal O}_X=0$ for $i > 0$ and $p_*{\cal O}_X\cong {\cal O}_{\PP 2}$
hence, by \cite{refHt}  III.8.3, $R^ip_* p^*\Omega (k+1)\cong
R^ip_*({\cal O}_X\otimes p^*\Omega (k+1))\cong R^ip_*{\cal O}_X
\otimes p^*\Omega (k+1)=0$ for $i > 0$ and $p_*p^*\Omega (k+1)
\cong \Omega (k+1)$, and so by \cite{refHt} ex.III.8.1 $H^i(p^*
\Omega (k+1)) \cong H^i(p_*p^*\Omega (k+1))= H^i(\Omega (k+1))$ for
all $i \geq 0$.

\begin{lem}\label{postlem}
If $Z$ has good postulation in degrees $k$ and $k+1$, then
$\mu_k$ is injective, resp. surjective, if and only if $\rho_k$ is
injective, resp. surjective.
Moreover, if $h^1({\cal I}_Z(k))= 0$, then
$$\hbox{exp-dim cok}\,\mu_k=
\hbox{exp-dim cok}\,\rho_k=\hbox{max}(0,\, 2l(Z)- k(k+2)),\hbox{ and}$$
$$\hbox{cok}\, \mu_k= \hbox{cok}\, \rho_k=H^1(\Omega (k+1)\otimes
{\cal I}_Z).$$
\end{lem}

\noindent{\bf Proof.} If $h^0({\cal I}_Z(k))=0$, $\mu_k$ is injective,
that is, $H^0(\Omega (k+1)\otimes {\cal I}_Z)=0$, so also $\rho_k$
is injective. If $h^1({\cal I}_Z(k))= 0$, we have that $\hbox{cok}\,\mu_k =
H^1(\Omega (k+1)\otimes {\cal I}_Z)= \hbox{cok}\,\rho_k$, and
$\hbox{ker}\,\mu_k
=H^0(\Omega (k+1)\otimes {\cal I}_Z)= \hbox{ker}\,\rho_k$, so that the
difference between the dimension of the domain and the
dimension of the codomain is the same: $h^0({\cal I}_Z(k))
h^0({\cal O}_{\PP 2}(1)) - h^0({\cal I}_Z(k+1))=3({k+2 \choose 2}-l(Z))-
({k+3 \choose 2}-l(Z))=k(k+2)-2l(Z)=h^0(\Omega (k+1))- h^0(\Omega
(k+1)\vert_{Z})$.
\prfend

\begin{lem}\label{equivlem}
The following conditions are equivalent:
\begin{description}
\item[(i)] $\rho_k$ is injective, resp. surjective;
\item[(ii)] $\psi_k$ is injective, resp. surjective;
\item[(iii)] $\eta_k$ is injective, resp. surjective.
\end{description}
\end{lem}

\noindent {\bf Proof.} $i) \Leftrightarrow ii)$:
one has $\pi_*({\cal E}_t)\cong \Omega(t)$,
and (see for example \cite{refId2}, 2.1)
$\pi_*({\cal E}_t\vert _T)\cong \Omega(t)\vert_Z$,
$\pi_*({\cal E}_t\otimes {\cal I}_T)\cong \Omega(t)\otimes {\cal I}_Z$;
hence $H^0({\cal E}_t)\cong H^0(\Omega(t))$, $H^0({\cal E}_t\otimes {\cal I}_T)
\cong H^0(\Omega(t)\otimes {\cal I}_Z)$, $H^0({\cal E}_t\vert
_T)\cong H^0(\Omega (t)\vert _Z).$

$i) \Leftrightarrow iii)$:
One has  $p_*{\cal O}_X \cong {\cal O}_{\PP2}$; by Prop. 2.3 of \cite{refAH},
one has also: $p_*{\cal O}_{\tilde Z}\cong {\cal O}_Z $;
hence it follows (see for example the proof of Lemma 2.3 in \cite{refId2},
taking into account that $p^{-1}(Z)={\tilde Z}$) that $p_*{\cal
I}_{\tilde Z}\cong {\cal I}_Z $.
By the projection formula we get $p_* (p^*\Omega(k+1)) \cong \Omega(k+1)$,
so $p_* (p^*\Omega(k+1)\vert_{{\tilde Z}} ) \cong \Omega(k+1)\vert_ {Z}$
and $p_* (p^*\Omega(k+1)\otimes {\cal I}_{\tilde Z}) \cong
\Omega(k+1)\otimes {\cal I}_Z$.
Hence the dimensions of the first three vector spaces in $(4*)$ and in
$ (2*)$ are the same, so we conclude that $\rho_k$ is of maximal rank
if and only if $\eta_k$ is.
\prfend

\section{Superfluous conditions for the cotangent bundle}\label{excondomega}

Now we are interested in studying the behaviour of the restriction of $\Omega
(k+1)$ to a curve in $\PP 2$. This will help us in the study of $\rho_k$ and
hence (see Section \ref{setup}) of $\mu_k$. In what follows $C$ will be a curve
of degree $d$ in $\PP 2$.

\begin{defn}\label{betadef}\rm  We denote by
$$\beta =\beta _{C,Z,k}: H^0(\Omega (k+1)\vert_C)\to H^0(\Omega
(k+1)\vert _{C\cap Z})$$
the restriction map. We also set
$$\gamma(C,Z,k):=\hbox{exp-dim cok}\,\beta_{C,Z,k}= max\{0,\; 2\,l(Z \cap
C)-h^0(\Omega (k+1)\vert_C)\}
.$$
\rm If $m(C)_{P_i}=r_i\leq m_i+1$, by Lemma \ref{lengthlem} $\; l(Z \cap C)=\sum
(r_im_i-{r_i\choose 2})$. So by Remark \ref{OmegaDimRem}, we find for $k+2\geq
d$ and $r_i\leq m_i+1$
$$\; \gamma(C,Z,k)=max\{0,\; 2\sum (r_im_i-{r_i\choose 2})-d(2k+2-d)\}.$$
\end{defn}

\begin{prop}\label{betasigmaprop}
Assume $h^1({\cal I}_Z(k))= 0$.
If $C\subset \PP2$ is a curve of degree $d \leq k+2$, then
$$\hbox{dim cok}\, \mu_k \geq  \hbox{dim cok}\, \beta_{C,Z,k}.$$
In particular, if there exists a (not necessarily integral) curve
$C$ of degree $d \leq k+2$ such that
$$\hbox{dim cok}\, \beta_{C,Z,k} > \hbox{exp-dim cok}\, \mu_k,$$
then  $\mu_k$ is not of maximal rank.
\end{prop}

\noindent{\em Proof}. Since $h^1({\cal I}_Z(k))= 0$, $Z$ has good postulation
in degree $k$ and $k+1$, so, by Lemma \ref{postlem},
$\hbox{dim cok}\, \mu_k = \hbox{dim cok}\, \rho _k = h^1({\cal
I}_Z\otimes \Omega (k+1))$.

Now set $t:=k+1$ and consider the commutative diagram:

$$\matrix {{} & {} & 0 & {} & 0 & {} & 0 & {} & {} \cr
{} & {} & \downarrow & {} & \downarrow & {} & \downarrow & {} & {} \cr
0 & \to & {\cal I}_{\hbox{res}_CZ}\otimes \Omega (t-d) & \to & {\cal
I}_Z\otimes
\Omega (t) & \to & {\cal I}_{Z\cap C,C}\otimes\Omega (t) & \to & 0 \cr
{} & {} & \downarrow & {} & \downarrow & {} & \downarrow & {} & {} \cr
0 & \to & \Omega (t-d) & \to & \Omega (t) & \to & {\cal O}_C\otimes
\Omega (t) & \to & 0 \cr
{} & {} & \downarrow & {} & \downarrow & {} & \downarrow & {} & {} \cr
0 & \to & {\cal O}_{\hbox{res}_CZ}\otimes \Omega (t-d) & \to & {\cal
O}_Z\otimes
\Omega (t) & \to & {\cal O}_{Z\cap C}\otimes \Omega (t) & \to & 0 \cr
{} & {} & \downarrow & {} & \downarrow & {} & \downarrow & {} & {} \cr
{} & {} & 0 & {} & 0 & {} & 0 & {} & {} \cr
}$$
Taking cohomology we get:
{\scriptsize
$$\matrix {{} & {} & 0 & {} & 0 & {} & 0 & {} & {} \cr
{} & {} & \downarrow & {} & \downarrow & {} & \downarrow & {} & {} \cr
0 & \to & H^0({\cal I}_{\hbox{res}_CZ}\otimes \Omega (t-d)) & \to &
H^0({\cal I}_Z\otimes \Omega (t)) & \to & H^0({\cal I}_{Z\cap
C,C}\otimes\Omega (t))
& \buildrel {\epsilon} \over \to & \cr
{} & {} & \downarrow & {} & \downarrow & {} & \downarrow & {} & {} \cr
0 & \to & H^0(\Omega (t-d)) & \to & H^0(\Omega (t))
& \to & H^0( \Omega (t)\vert_C) & \to & 0 \cr
{} & {} & \downarrow & {} & \downarrow\raise3pt\hbox
to0in{$\scriptstyle\rho_k$\hss} & {}
& \downarrow\raise2pt\hbox to0in{$\scriptstyle\beta_{C,Z,k}$\hss} & {} & {} \cr
0 & \to & H^0({\cal O}_{\hbox{res}_CZ}^{\oplus 2}) & \to & H^0({\cal
O}_Z^{\oplus 2})
& \to & H^0({\cal O}_{Z\cap C}^{\oplus 2}) & \to & 0 \cr
{} & {} & \downarrow & {} & \downarrow & {} & \downarrow & {} & {} \cr
{} & \buildrel {\epsilon} \over \to & H^1({\cal
I}_{\hbox{res}_CZ}\otimes \Omega (t-d))
& \to & H^1({\cal I}_Z\otimes \Omega (t)) & \to & H^1({\cal I}_{Z\cap
C,C}\otimes\Omega (t)) & \to & 0 \cr
{} & {} & \downarrow & {} & \downarrow & {} & \downarrow & {} & {} \cr
{} & {} & 0 & {} & 0 & {} & 0 & {} & {} \cr
}$$
}
where $h^1(\Omega (t)\vert_C)=0$ and $h^2({\cal I}_{\hbox{res}_CZ}\otimes
\Omega (t-d))=0$ since $h^2(\Omega (t-d))=0$ (this because $d \leq t+1$).
One has that $\hbox{dim cok}\, \mu_k = \hbox{dim cok}\, \rho _k =
h^1({\cal I}_Z\otimes
\Omega (t)) \geq h^1({\cal I}_{Z\cap C,C}\otimes\Omega (t))=
\hbox{dim cok}\, \beta_{C,Z,k}$.
\prfend

\par \medskip  As a consequence we have a first criterion to find schemes $Z$
for which $\mu_k$ fails to have maximal rank:
\begin{cor}\label{betasigmacor}
Assume $h^1({\cal I}_Z(k))= 0$. If there exists a (not necessarily integral)
curve $C$ of degree $d \leq k+2$ such that
\par \smallskip $\gamma(C,Z,k)>0$   if $\mu_k$ is exp-onto, or
\par \smallskip $\gamma(C,Z,k)>2l(Z)-k(k+2)$ if $\mu_k$ is exp-inj,
\b \smallskip  then  $\mu_k$ is not of maximal rank.
\end{cor}

\noindent{\em Proof}. This follows directly by Proposition \ref{betasigmaprop},
since $\hbox{dim cok}\, \beta_{C,Z,k} \geq \gamma(C,Z,k)$.
\prfend

\begin{example}\label{betasigmaex1} \rm
Let $Z=Z(3,2,1,1)$; then $\mu_4(Z)$ does not have maximal rank.
To see this directly, let $C$ be the line through $P_1$ and
$P_2$. Then $Z$ has good postulation (see \cite{refBHProc} or
\cite{refanti} for calculating the Hilbert function), $l(Z)=11$,
$h^0({\cal I}_Z(3))=0$, $h^0({\cal I}_Z(4))=4$, $h^0({\cal I}_Z(5))=10$, and
$\alpha=\tau=4$. Thus $I(Z)$ is generated in degrees at most 5, but since
$C$ is in the base locus of $H^0({\cal I}_Z(4))$ but the zero locus of
the whole ideal is just $Z$, there must be a generator of degree 5,
so the map $\mu _4(Z)$ is not surjective, and since it is exp-onto it is hence
not of
maximal rank (see \ref{CastelMumRegRem}).
\par Alternatively, note that
Corollary \ref{betasigmacor} applies: $\mu_4$ is exp-onto but
$\gamma(C,Z,4)>0$ since, using the fact
$\Omega (5)\vert_C \cong {\cal O}_C(3)\oplus {\cal O}_C(4)$,
we see $h^0(\Omega (5)\vert_C)=9<10=2\, l(Z\cap C)$. In other words, $Z\cap C$
imposes one superfluous condition to the sections of $\Omega (5)\vert_C$.
But a point of $\PP 2$ imposes 2 condition to a rank 2 bundle; so if we wish to
understand what's going on geometrically we have to move to $\PPP (\Omega)$.
Here
(cf. Lemma \ref{equivlem}) we have to
check the dimension of the space of global sections of
${\cal E}_5={\cal O}_{\PPP (\Omega)}(1)\otimes \pi^* {\cal O}_{\PP 2}(5)$
vanishing on the 1-dimensional scheme ${T}= \pi^{-1}(Z)$, or,
equivalently, on the 0-dimensional scheme
${T}':=({T} \cap \PPP (E))\cup ({T}\cap \PPP (G))$,
where $E\oplus G$ is a local trivialization of $\Omega$.
(In fact, since ${\cal E}_5$ is ${\cal O}_{\PP 1}(1)$ on the fibers,
the inverse
image of a point $\pi^{-1}(P)$ can be replaced by two generic points in
the fiber. For the non reduced case, and for further details, see
\cite {refId3}, \cite{refId2}, \cite{refGiId}.)
Hence we are looking at the postulation with respect to the invertible sheaf
${\cal E}_5$ of the 0-dimensional scheme ${T}'$ in $\PPP (\Omega)$;
since $\Omega (5)\vert_C \cong {\cal O}_C(3)\oplus {\cal O}_C(4)$,
we have a curve $D:=\PPP ({\cal O}_C(3))\subset \PPP
(\Omega)\vert_C\subset \PPP (\Omega)$
and ${T}'\cap D$ has length 5, while  ${\cal E}_5\vert_D\cong
{\cal O}_D(3)$ (see \cite{refHt} V.2.6).
It is now clear that it is possible to find a subscheme of $T'$ of length
$l({T}')-1$ imposing the same conditions as ${T}'$ to
${\cal E}_5$,
that is, ${T}'$ does not postulate well with respect to ${\cal E}_5$.
\end{example}

\begin{example}\label{betasigmaex2} \rm
Another similar example is given by
$Z=Z(4,3,3,3,2)$; here $\mu_7(Z)$ is again exp-onto and fails to have maximal
rank.
To see this, let $C$ be the conic through the 5 points $P_i$.
Again $Z$ has good postulation (\cite{refBHProc}, \cite{refanti}),
and we have $l(Z)=31$, $h^0({\cal I}_Z(6))=0$,
$h^0({\cal I}_Z(7))=5$, $h^0({\cal I}_Z(8))=14$, $\alpha=\tau=7$.
As before, $C$ is in the base locus of $I(Z)_7$, so
while the map $\mu_7$ is exp-onto it is not surjective, hence does
not have maximal rank.
Alternatively, again Corollary \ref{betasigmacor} applies: $\gamma(C,Z,7)>0$.
In more detail, $Z$ does not postulate well
with respect to $\Omega (8)$
(i.e., the number of sections of $\Omega (8)$ vanishing on $Z$
is greater than the length of $Z$ would lead us to expect),
since $h^0(\Omega (8)_{\vert C})=28<2\cdot15=2l(Z\cap C)$.
In other words, $Z\cap C$ imposes 2 superfluous conditions
on the sections of $\Omega (8)\vert_C\cong {\cal O}_{\PP 1}(13)^{\oplus 2}$
(here we have used the fact that $\Omega\vert_C\cong {\cal O}_{\PP
1}(-3)^{\oplus 2}$).
If we work in $\PPP (\Omega)$ (cf. Lemma \ref{equivlem}),
we have to consider the postulation with respect to ${\cal E}_8$ of
the 1-dimensional
scheme ${T}= \pi^{-1}(Z)$, or, as in the previous example,
of the 0-dimensional scheme ${T}':=({T}\cap \PPP
(E))\cup ({T}\cap \PPP (G))$,
$E\oplus G$ again being a local trivialization of $\Omega$.
Since $\Omega (8)\vert_C\cong {\cal O}_{\PP 1}(13)^{\oplus 2}$, we
have two curves,
$D_1$ and $D_2$,  both contained in $\PPP(\Omega)\vert_C$, where
${T}'\cap D_i$ has length 15, while  ${\cal
E}_8\vert_{D_i}\cong {\cal O}_{\PP 1}(13)$.
It is hence possible to find a subscheme of ${T}'$ of length
$l({T}')-2$
imposing the same conditions as $T'$ on ${\cal E}_8$; i.e., ${T}'$
does not postulate well with respect to ${\cal E}_8$.
These two superfluous conditions give a contribution of 2 to the cokernel.
\end{example}

 \begin{example}\label{3triple} \rm
The map $\mu_{5}$ for $Z=3P_1+3P_2+3P_3$ fails to have maximal rank; $Z$
postulates well (\cite{refHi1}), $l(Z)=18$,
$h^0({\cal I}_Z(4))=0$, $h^0({\cal I}_Z(5))=3$, $h^0({\cal I}_Z(6))=10$,
$\alpha=\tau=5$, and $\mu_5$ is exp-inj. But $\mu_5$ is not injective; actually,
if $L_{ij}$ is the line through $P_i$ and $P_j$, and $C$ is the union of
$L_{12}$, $L_{13}$ and $L_{23}$, the cubic $C$ is a fixed component for $\vert 
I(Z)_5\vert$, hence the three generators of $ I(Z)_5$ are of the form $CF_i$,
$i=1,2,3$, where $F_1$, $F_2$ and $F_3$ are the three conics which generate
$I(P_1+ P_2+ P_3)$. Since $h^0({\cal I}_{P_1+ P_2+ P_3}(3))=7$, the dimension
of the image of $\mu_5$ is also 7, i.e. $I(Z)$ needs 3 generators in degree 6,
not just one.
\par Alternatively, note that Corollary \ref{betasigmacor} applies:
$\gamma(C,Z,5)=3>2l(Z)-5(5+2)=1$; here $C$ is a triangle, hence reducible with
aritmethic genus 1.
What happens here is that $\Omega (6)\vert_{L_{ij}} \cong {\cal
O}_{L_{ij}}(4)\oplus {\cal O}_{L_{ij}}(5)$, so that $Z\cap L_{ij}$ imposes one
superfluous condition to the sections of $\Omega (5)\vert_{L_{ij}}$ for each
one of the three lines $L_{ij}$; one of these superfluous conditions wouldn't
bother the rank maximality of $\mu_{5}$, which is expected to be injective with
a 1-dimensional cokernel; the other two conditions give a contribution of 2 to
the cokernel.
If we reinterpret the situation in $\PPP (\Omega)$, we have to consider the
postulation with respect to ${\cal E}_6$ of
a certain 0-dimensional scheme $T'$, analogously to what happens in the previous
examples; here there are three curves $D_{ij}:=\PPP ({\cal O}_{L_{ij}}(4))$ such
that ${\cal E}_6\vert_{D_{ij}}\cong {\cal O}_{D_{ij}}(4)$, while ${T}'\cap
D_{ij}$ has length 6; ${T}'$ does not postulate well with respect to ${\cal
E}_6$. Notice anyway that the reducible curve $D_{12}\cup D_{13}\cup D_{23}$
causing troubles is now the union of three disjoint smooth rational curves,
since a point $P$ in $D_{ij}$ is the point in $\PPP (\Omega \vert _{L_{ij}})$
representing the tangent direction of $L_{ij}$ at $P$.
\end{example}

These first three examples are easy to treat by taking into account the
occurrence of fixed components. The next example (as well as example
\ref{betasigmaex4}) shows that this is not always the case.

\begin{example}\label{betasigmaex3} \rm
If $Z= 9P_1+\cdots+9P_7$, the map $\mu_{24}$ fails
to have maximal rank. Again $Z$ postulates well (\cite{refIGP},
\cite{refBHProc},
\cite{refanti}), $l(Z)=315$,
$h^0({\cal I}_Z(23))=0$, $h^0({\cal I}_Z(24))=10$, $h^0({\cal I}_Z(25))=36$,
$\alpha=\tau=24$, and $\mu_{24}$ is exp-inj.
But $\vert I(Z)_{24} \vert$ is fixed component free and $\mu_{24}$ is not
injective; if it were, $\hbox{dim Im}\mu_{24}$
would be 30,
but in fact it is $29$  (\cite{refIGP}).
Once more, Corollary \ref{betasigmacor} applies:
$\gamma(C,Z,24)>2l(Z)-24(24+2)$,
with $C:= \sum C_i$, where $C_i$ is a cubic with $m(C_i)_{P_j}=1$ for
$i\neq j$, 2 for $i=j$.
Again, $C$ is not irreducible.
What happens here is that the superfluous conditions imposed by
$Z\cap C$ on $\Omega (25)\vert_C$
are more than the expected dimension for the cokernel of $\rho_{24}$, since
$2l(Z \cap C)- h^0(\Omega (25)\vert_C)=616-609=7>2l(Z)-24(24+2)=6$.
Thus $\rho_{24}$, and hence $\mu_{24}$, are not injective by
Corollary \ref{betasigmacor}.
(Notice that taking into account just one of the curves $C_i$ is not enough:
in fact, $h^0(\Omega (25)\vert_{C_i}) =141 < 2\, l(Z \cap C_i)=142$;
but this only says that $\hbox{dim cok} \rho_{24} \geq 1$.)
\end{example}

\section{Superfluous conditions for the pullback of the cotangent
bundle}\label{excondpullback}

In examples \ref{betasigmaex1}, \ref{betasigmaex2}, \ref{3triple} and
\ref{betasigmaex3},
failure of $\mu_k(Z)$ to have maximal rank was related to $Z$ imposing too many
conditions on the global sections of $\Omega (k+1)\vert_C$, and we checked it
just by a dimension count, i.e. the expected dimension $\gamma(C,Z,k)$ of the
cokernel of $\beta_{C,Z,k}$ was too big. But $\Omega (k+1)\vert_C$ is a {\em
rank two} vector bundle, so it can happen that the dimension of the cokernel is
bigger than its expected dimension. This of course cannot occur
with a rank one bundle on $\PP1$, since if $A$ is a 0-dimensional scheme on
$\PP1$, the cokernel of the restriction map $H^0({\cal O}_{\PP1}(t))\to
H^0({\cal O}_A)$ always has the expected dimension.
\par Instead if we consider for
example the restriction map $H^0({\cal O}_{\PP1}\oplus {\cal O}_{\PP1}(2))\to
H^0({\cal O}_A^{\oplus 2})$ where $A$ is the union of two points, then the
expected
dimension of the cokernel is $0$ but the actual dimension is 1; $A$ imposes 1
condition too many on $H^0({\cal O}_\PP1)$. This is possible because the
splitting gap of
${\cal O}_{\PP1}\oplus {\cal O}_{\PP1}(2)$ is 2. In the previous examples this
behaviour did not arise: in examples \ref{betasigmaex1} and
\ref{betasigmaex2}, $C$ is a
line or a smooth conic with splitting gap 1, respectively 0; in example
\ref{betasigmaex3},
$C_i$ is a singular cubic, so we don't look at $\Omega (k+1)\vert_C{_i}$, but
we can look at the splitting of the pull-back of $\Omega(k+1)$ on $\tilde C_i$,
and we find that the splitting gap is 1.
\par The forthcoming example \ref{betasigmaex4}, instead, illustrates a
situation such that
$\mu_k(Z)$ is exp-onto, but there exists a curve $C$ with  splitting gap 2, and
$\hbox{dim cok}\,\beta_{C,Z,k}> 0$, so $\mu_k(Z)$ is not onto although
$\hbox{exp-dim  cok}\,\beta_{C,Z,k}=\gamma(C,Z,k)=0$.
So it seems evident that, if we want to formulate a conjecture about the rank
maximality of $\mu_k(Z)$, it is necessary to take into consideration the
splitting type, and to consider the real cokernel of the maps $\,\beta_{C,Z,k}$;
this is what we are going to do next.

\begin{defn}\label{delta0} \rm Let $C$ be a curve of degree $d$ in $\PP 2$, such
that its strict
transform $\tilde C= dL-\sum r_iE_i$ is smooth and rational in the surface
$X$ obtained by blowing up the points $P_i$.  Given a positive integer $k$ and
taking cohomology of the exact sequence
$0 \to p^*\Omega (k+1)\vert_{\tilde C}\otimes {\cal I}_{\tilde C \cap
\tilde Z}\to
p^*\Omega (k+1)\vert_{\tilde C}\to p^*\Omega (k+1)\vert_{\tilde C
\cap \tilde Z}\to 0$, where $\tilde Z=\sum m_iE_i$, we get the restriction map

$$\theta = \theta_{C,Z,k}:  H^0(p^*\Omega (k+1)\vert_{\tilde C})
\to H^0(p^*\Omega (k+1)\vert_{\tilde C \cap \tilde Z}).$$
In order to measure the superabundance of conditions imposed by $\tilde C \cap
\tilde Z$ on the sections of $p^*\Omega (k+1)\vert_{\tilde C}$ we also set

$$\delta_0(C,Z,k)=\hbox{dim cok}\,\theta _{C,Z,k}.$$

Writing $a$ and $b$ for $a_C$ and $b_C$, we have
$p^*\Omega (k+1)\vert_{\tilde C}\cong  {\cal O}_{\tilde C}(-a+dk)\oplus
{\cal O}_{\tilde C}(-b+dk)$. Moreover, $\tilde C \cdot \tilde
Z =\sum r_i m_i$.

Since $b\leq d$ and by assumption $k\geq 1$, we have $dk-b\geq0$ so that
$h^1(p^*\Omega (k+1)\vert_{\tilde C})=0$.
Hence $\delta_0(C,Z,k)=h^1(p^*\Omega (k+1)\vert_{\tilde C}\otimes {\cal
I}_{\tilde C \cap \tilde Z})
=h^1({\cal O}_{\PP 1}(-a+dk-\sum r_i m_i)\oplus {\cal O}_{\PP
1}(-b+dk-\sum r_i m_i))
=\hbox{max}(0,\; l(\tilde Z\cap \tilde C)-h^0({\cal O}_{\tilde C}(-a+dk))
+ \hbox{max}(0,\; l(\tilde Z\cap \tilde C)-h^0({\cal O}_{\tilde C}(-b+dk))$, so
that finally
$$\delta_0(C,Z,k)=\hbox{max}(0,\; \sum r_im_i-dk+a-1) + \hbox{max}(0,\; \sum
r_im_i-dk+b-1).$$
\end{defn}

\medskip In certain cases, $\delta_0$ is nothing more than $\gamma$:

\begin{thm}\label{betarholem}
Let $C\subset \PP 2$ be a curve whose strict transform $\tilde
C=dL-\sum r_iE_i$ is smooth
and rational in $X$, and assume $d\leq k+2$ and $r_i-1 \leq m_i$ for all $i$.
Then
$$\hbox{cok}\, \beta_{C,Z,k}\cong \hbox{cok}\,\theta_{C,Z,k}$$
hence $\delta_0(C,Z,k)\geq \gamma(C,Z,k)$, with equality
if and only if $\hbox{cok}\, \beta_{C,Z,k}$ has its expected dimension
(this occurs, for example, if $\sum r_im_i-dk+a-1\geq 0$).
\end{thm}

\noindent{\em Proof}. The maps $\theta=\theta_{C,Z,k}$ and
$\bar{\theta} =\bar{\theta}_{C,Z,k}: H^0(p_*p^*\Omega (k+1)\vert_{\tilde C})
\to H^0(p_*p^*\Omega (k+1)\vert_{\tilde C \cap \tilde Z})$
are the maps on cohomology coming from the exact sequences
$0 \to p^*\Omega (k+1)\vert_{\tilde C}\otimes {\cal I}_{\tilde C \cap
\tilde Z}\to
p^*\Omega (k+1)\vert_{\tilde C}\to p^*\Omega (k+1)\vert_{\tilde C
\cap \tilde Z}\to 0$
and its pushforward by $p_*$, and it is clear that $\hbox{cok}\,
\bar{\theta}\cong \hbox{cok}\,\theta.$

We also have an exact sequence $0 \to {\cal O} _C\to  p_*{\cal
O}_{\tilde C}\to {\cal S} \to 0$,
where ${\cal S} = \oplus_{P\in \hbox{Sing}(C)}{\tilde {\cal
O}}_P/{{\cal O}}_P$ and
${\tilde {\cal O}}_P$ denotes the integral closure of ${\cal O}_P$.
Letting $\delta _P$ be the length $l({\tilde {\cal O}}_P/{ {\cal O}}_P)$,
one has (by \cite{refHt}, Ex. IV.1.8 and Cor. V.3.7)
$p_a(C)=p_a(\tilde C)+\sum_{P\in \hbox{Sing}(C)}\delta _P$.
But $0=p_a(\tilde C)= {d-1 \choose 2}-\sum {r_i \choose 2}$ and
$p_a(C) ={d-1 \choose 2}$, so
$\sum_{P\in\hbox{Sing}(C)}\delta _P=  \sum {r_i \choose 2}$, hence
$l({\cal S})=\sum {r_i \choose 2}$.
\b There is a natural map
$0 \to {\cal O}_{ C \cap Z}\to  p_*{\cal O}_{\tilde C \cap \tilde Z}$; let us
denote the cokernel by ${\cal S'}$. Since $r_i\leq m_i+1$ for all $i$, by lemma
\ref{lengthlem} we have $l({\cal S'})=l(\tilde C \cap
\tilde Z)-l(C\cap Z)= \sum r_im_i- \sum (r_im_i-{r_i\choose 2})=\sum {r_i
\choose 2}$. Now consider the diagram

$$\matrix {
0 & \to & {\cal O} _C & \to & p_*{\cal
O}_{\tilde C} &  \to & {\cal S} & \to & 0 \cr
{} & {} & \downarrow\raise3pt\hbox to0in{$\scriptstyle $\hss} &
{} & \downarrow\raise3pt\hbox to0in{$\scriptstyle $\hss}
& {} & \raise3pt\hbox to0in{$\scriptstyle $\hss} & {} & {} \cr
0 & \to & {\cal O}_{ C \cap Z} & \to &  p_*{\cal O}_{\tilde C \cap \tilde Z}&
\to
& {\cal S'} & \to & 0 \cr
{} & {} & \downarrow\raise3pt\hbox to0in{$\scriptstyle $\hss} &
{} & \downarrow\raise3pt\hbox to0in{$\scriptstyle $\hss}
& {} &  \raise3pt\hbox to0in{$\scriptstyle $\hss} & {} & {} \cr
{} & {} & 0 & {} & R^1p_*{\cal
I}_{\tilde C\cap \tilde Z, \tilde C} &  {} & {} & {} & {} \cr
}$$

There is a map ${\cal S}\to {\cal S'}$ making the diagram commute, and it
has to be surjective since $R^1p_*{\cal I}_{\tilde C\cap \tilde Z, \tilde C}=0$
by \cite{refHt} III.11.2. Hence it is a surjective map between sheaves
supported at points and of the same lenght, so we conclude ${\cal S'}\cong
{\cal S}$, which gives us the exact sequence $0 \to {\cal O}_{ C \cap Z}\to
p_*{\cal O}_{\tilde C \cap \tilde Z}\to {\cal S'}\to 0$.

Tensoring this and the exact sequence at the beginning of the proof
by $\Omega (k+1)$, taking into account that
$p_*{\cal O}_{\tilde C}\otimes \Omega (k+1)\cong p_*p^*\Omega
(k+1)\vert_{\tilde C}\,$
and $\, p_*{\cal O}_{\tilde C \cap \tilde Z}\otimes \Omega (k+1)\cong
p_*p^*\Omega (k+1)\vert_{\tilde C \cap \tilde Z}$
(cf. the projection formula, \cite{refHt} III.8.3), recalling that $H^1(\Omega
(k+1)_{\vert C})=0$ for $k+2\geq d$ (see Remark \ref {OmegaDimRem}), and
finally writing
$\beta=\beta _{C,Z,k}$, we get
$$\matrix {
0 & \to & H^0(\Omega (k+1)_{\vert C}) & \to & H^0(p_*p^*\Omega
(k+1)\vert_{\tilde C}) &  \to & H^0({\cal S}^{\oplus 2}) & \to & 0 \cr
{} & {} & \downarrow\raise3pt\hbox to0in{$\scriptstyle\beta$\hss} &
{} & \downarrow\raise3pt\hbox to0in{$\scriptstyle\bar{\theta}$\hss}
& {} & \downarrow\raise3pt\hbox to0in{$\scriptstyle\cong$\hss} & {} & {} \cr
0 & \to & H^0(\Omega (k+1)\vert_{C\cap Z}) & \to & H^0(p_*p^*\Omega
(k+1)\vert_{\tilde C \cap \tilde Z}) & \to
& H^0({\cal S}^{\oplus 2}) & \to & 0 \cr
}$$
The snake lemma now gives  $\hbox{cok}\,\bar{\theta}\cong \hbox{cok}\,\beta$.

The inequality $\delta_0(C,Z,k)\geq \gamma(C,Z,k)$ is now clear,
since $\delta_0$ is the dimension of $\hbox{cok}\,\theta_{C,Z,k}$,
while $\gamma$ is merely the expected dimension of
$\hbox{cok}\,\beta_{C,Z,k}$. For the rest, assuming
$\sum r_im_i-dk+b-1 \geq \sum r_im_i-dk+a-1\geq 0$ and using
$h^0(p^*\Omega (k+1)\vert_{\tilde C})=h^0(p_*p^*\Omega (k+1)\vert_{\tilde C})
=h^0(\Omega (k+1)\vert _C+2\sum{r_i \choose 2}$,
we have $\delta_0(C,Z,k)= 2l(\tilde Z\cap \tilde C)-h^0(p^*\Omega
(k+1)\vert_{\tilde C})=
2(l(\tilde Z\cap \tilde C)-\sum{r_i \choose 2})-
h^0(\Omega (k+1)\vert_C=2\,l(Z \cap C)-h^0(\Omega (k+1)\vert_C)=
\gamma (C,Z,k)$.
\prfend

\begin{cor}\label{betarhocor}
Assume $h^1({\cal I}_Z(k))= 0$ and moreover that there exists an integral curve
$C\subset \PP 2$ such that $\tilde C=dL-\sum r_iE_i$ is smooth and rational in
$X$ with $d\leq k+2$, $r_i-1 \leq m_i$ for all $i$ and
$\delta_0(C,Z,k)>\hbox{exp-dim cok}\, \mu_k(Z)$.
Then $\mu_k(Z)$  is not of maximal rank.
\end{cor}

\noindent{\em Proof}. We have
$\hbox{exp-dim cok}\, \mu_k<\delta _0(C,Z,k)
=\hbox{dim cok}\, \theta _{C,Z,k}= \hbox{dim cok}\,\beta_{C,Z,k}$
and we conclude by Proposition \ref{betasigmaprop}.
\prfend

We now show how to use this last result.
A significant difference here with the three previous examples
is that the splitting gap for (any irreducible component of) $C$ was 0 or 1
previously; in Example \ref{betasigmaex4} it is 2.

\begin{example}\label{betasigmaex4} \rm
Let $Z= 4P_1+\cdots+4P_7+P_8$; then $\mu_{11}(Z)$ fails to
have maximal rank (see \cite{refFHH}).
Note that  $Z$ has good postulation and
$I(Z)_{11}$ is fixed component free (apply \cite{refBHProc} or \cite{refanti}).
We have  $l(Z)=71$,  $h^0({\cal I}_Z(10))=0$, $h^0({\cal I}_Z(11))=7$,
$h^0({\cal I}_Z(12))=20$, $\alpha=\tau=11$, hence $\mu_{\alpha}$ is exp-onto.
The map $\mu_{11}(Z)$ is
not surjective.
This can be attributed to to the existence of a rational curve $C$ of
degree 8 with $r_i:=m(C)_{P_i}=3$ for $0\leq i \leq 7$, and $r_8=1$;
$\tilde C\subset X$ is a smooth rational curve of self-intersection
$\tilde C^2=0$.
This time we cannot read failure of maximal rank on the sections of $\Omega$,
since $h^0(\Omega (12)_{\vert C})=128=2\, l(Z\cap C)$; i.e., $\gamma(C,Z,11)=0$.
Instead, the splitting gap for $C$ is 2, since  (see the proof of
Lemma 12 of  \cite{refFHH})
$p^*(\Omega (1))\vert_{\tilde C} \cong {\cal O}_{\tilde C}(-3)\oplus
{\cal O}_{\tilde C}(-5)$,
hence $p^*(\Omega (12))\vert_{\tilde C} \cong {\cal O}_{\tilde C}(85)\oplus
{\cal O}_{\tilde C}(83)$.
The scheme $\tilde Z:= \sum m_iE_i$ intersects $\tilde C$ in a 0-dimensional
scheme of length $\sum r_im_i=85$, so $\tilde Z \cap \tilde C$ is too much for
${\cal O}_{\PP 1}(83)$ (and not enough for ${\cal O}_{\PP 1}(85)$);
that is, the cohomology
of the exact sequence
$0 \to p^*\Omega (12)\vert_{\tilde C}\otimes {\cal I}_{\tilde C \cap \tilde Z}
\to p^*\Omega (12)\vert_{\tilde C} \to p^*\Omega (12)\vert_{\tilde C
\cap \tilde Z} \to 0$ is

$$0\to H^0({\cal O}_{\tilde C}\oplus {\cal O}_{\tilde C}(-2))\to
H^0(p^*\Omega (12)\vert_{\tilde C})
\buildrel {\theta } \over \to H^0(p^*\Omega (12)\vert_{\tilde C \cap \tilde Z})
\to  H^1({\cal O}_{\tilde C}\oplus {\cal O}_{\tilde C}(-2)) \to 0$$

\b so $\theta $ is not of maximal rank. (This cannot happen if the
splitting gap is 0 or 1.)
Since $\delta_0(C,Z,11)=1>0$, we see by Corollary \ref{betarhocor} that
$\mu_{11}(Z)$ fails to have maximal rank.
\end{example}

\bigskip The previous examples might lead one to think that the curves $C$ that
need to be taken into consideration are the ones with $C^2\leq 0$. The
following example shows that this is not the case.

\begin{example}\label{autointex} \rm
Let $Z=15(P_1+\dots +P_4)+13(P_5+P_6)+9P_7+2(P_8+\dots+P_{11})$; then
$l(Z)=719$, $h^0({\cal I}_Z(37))=22$, $h^0({\cal I}_Z(38))=61$, so that
$\mu_{37}(Z)$ is exp-onto but in fact it does not have maximal rank; precisely,
$\hbox{dim cok}(\mu_{37})=1$; this has been computed with Macaulay 2
(\cite{refGS}).
\b Now consider a curve $C$ whose strict transform $\tilde C$ is an irreducible
curve in the linear system $\vert 34L- 14(E_1+\dots
+E_4)-12(E_5+E_6)-8E_7-2(E_8+\dots+E_{11})\vert $ (such a $C$ exists, since
$\tilde C=2D$ with $D$ a Cremona transform of a line, hence the linear system
above contains Cremona transforms of conics). One has $\tilde C^2=4$. The
splitting type for $\tilde C$ is $(14,20)$ (it is possible to compute it with
the script in \cite{refGHIpv}); then (\ref{delta0}) $\delta_0(C,Z,37)=1$. Here
too we cannot work in $\PP 2$; in fact,  $\gamma(C,Z,37)=0$.

\end{example}

\section{Two conjectures}\label{twoconj}

In each of our examples above, failure of $\mu_k(Z)$ to be surjective is
accompanied by $\delta_0(C,Z,k)>0$. This seems to be fairly general
behavior, which leads us to advance the following conjecture.

\begin{conj}\label{conj1}
Let $Z=\sum m_iP_i$ be a fat point scheme in $\PP 2$ (for general
points $P_i$),
with $h^1({\cal I}_{Z}(k))=0$. Say $\mu_k(Z)$ is exp-onto. Then
$\mu_k(Z)$ fails to be surjective if and only if there exists an
integral curve $C\subset \PP 2$ whose strict transform $\tilde C=
dL-\sum r_iE_i$
is smooth and rational in $X$, with $d\leq k+2$, $\; r_i \leq m_i+1$ and
$\delta_0(C,Z,k) >0$.
\end{conj}

\begin{remark}\label{rmconj1} \rm
In fact, in every example we have found for which $\mu_k(Z)$ fails to have
maximal rank, we have $h^0(kL-\sum m_iE_i-\tilde C)>0$, and hence $d\leq
k$.\end{remark}

\b The ``if" part of Conjecture \ref{conj1} is true, and is
Corollary \ref{betarhocor}. Here are some counterexamples to the ``if" part of
conjecture \ref{conj1} with $d> k+2$ and $ r_i > m_i+1$ for some $i$:

$Z=P_1$, $k=1$, $\tilde C=4L-3E_1-E_2-\cdots-E_8$;

$Z=P_1+\cdots+P_4$, $k=2$, $\tilde C=5L-3E_1-2(E_2+E_3+E_4)-(E_4+\cdots+E_8)$;

$Z=P_1+\cdots+P_7$, $k=3$, $\tilde C=8L-3(E_1+\cdots+E_7)-E_8$;

$Z=2P_1+2P_2+P_3+\cdots+P_7$, $k=4$, $\tilde
C=7L-4E_1-3E_2-2(E_3+\cdots+E_8)$.
\b The problem in each case is, in some sense, that $C$ is too big.

\bigskip In the case when $\mu_k(Z)$ is exp-inj the situation is more
complicated. We have already seen in Examples \ref{3triple}, \ref{betasigmaex3}
that the curve $C$ needs not be irreducible; the following example shows that it
can also be nonreduced.

\begin{example}\label{nonreducedex} \rm
Let $Z=60(P_1+\dots +P_8)$; then $l(Z)=14640$, $h^0({\cal I}_Z(169))=0$,
$h^0({\cal I}_Z(170))=66$, $h^0({\cal I}_Z(171))=238$, $\alpha=\tau=170$, so
that $\mu_{170}(Z)$ is exp-inj but in fact it does not have maximal rank (see
\cite{refIGP}); precisely, $\hbox{exp-dim cok}(\mu_{170})=40$, while the actual
dimension is $48$.
\b Let $C_j$ be a sextic with $r_{j,i}=m(C_j)_{P_i}=2$ for
$i\neq j$ and $r_{j,j}=m(C_j)_{P_j}=3$,  $\, j=1,\dots,8$. The splitting type
for $C_i$ is $(3,3)\,$ by \ref{splitlem}
; then (\ref{delta0}) $\delta_0(C_j,Z,170)=2\, \hbox{max}(0,\; 60\,\sum _i
r_{j,i}-6\cdot170+3-1)=4$. In order to take into account the contribution of
each $C_j$, we do as in Example \ref{betasigmaex3} and we consider  $C=\sum
C_j$, but this is still not enough since $8\cdot4<\hbox{exp-dim
cok}(\mu_{170})$.

\b So we go on: since $res_{C_j}Z=57P_j+58\,\sum _{i\neq j}P_i$, we find
$\delta_0(C_j,res_{C_j}Z,170-6)=2$. If we add up the contribution not only of
$Z$ but also of $res_{C_j}Z$ for all the $C_j$'s, we then find $\hbox{dim
cok}(\mu_{170})\geq 8(4+2)=48$. It is useless to go on, since
$\delta_0(C_j,res_{C_j}(res_{C_j}Z),170-12)=0$.
\b Notice that since the splitting type for $C_i$ is balanced, we can work
directly in $\PP 2$; it is easy to check that  $\gamma(2C,Z,170)=48$ so it is
enough to apply Proposition \ref{betasigmaprop}.

\end{example}

\par In this last example we have seen that it is enough to consider $\gamma$,
but this is not always the case for injectivity too. In fact, in the following
example bijectivity is expected, and $\gamma =0$, while $\delta_0=1$.

\begin{example}\label{delta0inj} \rm
Let $Z=11(P_1+\dots +P_7)+5P_8+2P_9$; then $l(Z)=480$, $h^0({\cal I}_Z(30))=16$,
$h^1({\cal I}_Z(30))=0$, $h^0({\cal I}_Z(31))=48$, so that $\mu_{30}(Z)$ is
exp-bijective but in fact it does not have maximal rank. To see this, it is
enough to apply Corollary \ref {betarhocor} with  $\tilde C =19L-7(E_1+\dots
+E_7)-4E_8-E_9$.  The splitting type for $C$ is $(8,11)$ (to compute it, use
\cite{refGHIpv}); then \ref{delta0} gives $\; \delta_0(C,Z,30)=\hbox{max}(0,\;
-9+8-1) + \hbox{max}(0,\; -9+11-1)=1$.

\b On the other hand, it is easy to check (\ref{betadef}) that $\;
\gamma(C,Z,30)=0$, so Corollary \ref{betasigmacor} is useless here.

\end{example}

These examples motivate the following definition:

\begin{defn}\label{delta} \rm Let $C\subset \PP 2$ be a degree $d$ curve, with
$m(C)_{P_i}=r_i$.
The $h$-iterated residual scheme of $Z=\sum m_iP_i$ with respect to
$C$ is defined inductively as follows:
$$res_{C,0}Z:=Z, \quad res_{C,h}Z:=res_C(res_{C,h-1}Z).$$
\b Notice that $res_{C,h}Z=\sum (m_i-hr_i)P_i$ if $m_i\geq hr_i$.

\b Assume now that the strict transform $\tilde C= dL-\sum r_iE_i$ is
smooth rational with $a:=a_C,\; b:=b_C$. Let $t-hd \geq 1$; we
define inductively the $h$-superabundance of $C$:

$$\delta _h(C,Z,t):= \delta _0(C,res_{C,h}Z,\, t-hd).$$

\b We finally set $$\delta (C,Z,t):= \sum_{h=0,...,[{t\over d}]} \delta
_h(C,Z,\, t).$$

Now let $F$ be as usual $F=tL-\sum m_iE_i$, and set $A _h(C,Z,\, t)=-F\cdot
\tilde C+a-1+h\,\tilde C^2$, $B _h(C,Z,\, t)=-F\cdot \tilde C+b-1+h\,\tilde
C^2$. Then $B _h(C,Z,\, t)\geq A _h(C,Z,\, t)$, and if $m_i\geq hr_i$, $t-hd
\geq 1$, we have:
\b $\delta_h(C,Z,t)=\hbox{max}(0,\; \sum r_i(m_i-hr_i)-d(t-hd)+a-1) +
\hbox{max}(0,\; \sum r_i(m_i-hr_i)-d(t-hd)+b-1)=\hbox{max}(0,\; A_h(C,Z,\, t))
+ \hbox{max}(0,\; B_h(C,Z,\, t))$.
\end{defn}

\a To understand better the connection between $\delta$ and $\gamma$, the
following proposition is helpful:

\begin{prop}\label{recursive} Let $C\subset \PP 2$ be a curve with $\tilde C=
dL-\sum r_iE_i$ smooth rational.
 Assume that $(p+1)r_i-1\leq m_i$, $t+2\geq d(p+1)$, and assume also that
$A_h(C,Z,\, t) \geq 0$ for $h=0,\ldots,p$. Then, denoting by $(p+1)C$ the
$p^{th}$
infinitesimal neighborhood of $C$ in $\PP 2$, one has:
$$\sum_{h=0,...,p} \delta _h(C,Z,\, t) =\gamma (\,(p+1)C,Z,t).$$
\end{prop}

\noindent{\em Proof}. First notice that, if $A_0(C,Z,\, k)\geq 0$, then using
adjunction formula $\delta_0(C,Z,k)=A_0(C,Z,\, k) + B_0(C,Z,k)=2\sum
(r_im_i-{r_i\choose 2})-d(2k+2-d)\geq 0$. Hence, if $k+2\geq d$ and $r_i\leq
m_i+1$, then
$\; \gamma(C,Z,k)= 2\sum (r_im_i-{r_i\choose 2})-d(2k+2-d)=2\,l(Z \cap
C)-h^0(\Omega (k+1)\vert_C)=\delta_0(C,Z,k)$ (see \ref {betadef},
\ref{betarholem}).
\par We have $res_{C,h}Z=\sum(m_i-hr_i)P_i$, since by assumption $hr_i\leq m_i$,
 for $0\leq h \leq p$.
\par So, since by assumption  $r_i-1\leq m_i-hr_i\,$, $t-hd+2\geq d$ and
$A_h(C,Z,\, t) \geq 0$  for $h=0,\ldots,p$, we have $\delta_h(C,Z,t)=\delta
_0(C,\sum(m_i-hr_i)P_i,\, t-hd)=\gamma (C,\sum(m_i-hr_i)P_i,\,
t-hd)=2\,l((\sum(m_i-hr_i)P_i) \cap C)-h^0(\Omega (t-hd+1)\vert_C)$  for
$h=0,\ldots,p$.
\par It is easy to check (see \ref {OmegaDimRem}, \ref {lengthlem} and use
$t-dh+2\geq d$ and $r_i-1\leq m_i-hr_i$ for $0\leq h \leq p$, and
$m((p+1)C)_{P_i}=(p+1)r_i \,$) that
\b $\sum_{h=0,...,p}h^0(\Omega(t-hd+1)_{\vert
C})=h^0(\Omega(t+1)_{\vert  (p+1)C})\;$,  and
\b $\sum_{h=0,...,p}l((\sum (m_i-hr_i)P_i)\cap C)=l((\sum m_iP_i)\cap
(p+1)C)$. So conclusion follows adding up.
\prfend

\par We are now ready to formulate a conjecture for the case where injectivity
is expected.

\begin{conj}\label{conj2}
Let $Z=\sum_i m_iP_i$ be a fat point scheme in $\PP 2$  (for general
points $P_i$), with $h^1({\cal I}_{Z}(k))=0$. Say $\mu_k(Z)$ is exp-inj. Then
$\mu_k(Z)$ fails to be injective if and only if there exists a curve
$C\subset \PP 2$ such that: $\tilde C= dL-\sum r_iE_i$ has $r_i\leq m_i+1$ and
$d\leq k+2$; $ C=\sum n_j C_j$, where each $C_j$ is integral with  $\tilde C_j$
smooth and rational in $X$ and $\tilde C_{j_1}\cdot
\tilde C_{j_2}=0$; and $\sum \delta (C_j,Z,k) >2l(Z)-k(k+2)$.
\end{conj}

\b The ``if" part of conjecture \ref{conj2} is true if for example $j=1$ and
$A_h(C,Z,\, t) \geq 0$ for $h=0,\ldots,n_1$ by Proposition \ref {recursive} and
Corollary \ref {betasigmacor}.

\bigskip \par Notice that all the results on the generation for fat point
schemes (see the introduction for a list of them) are consistent with
Conjectures  \ref {conj1} and \ref{conj2}.

\bigskip We end by proving that the SHGH Conjecture together with
Conjectures \ref{conj1} and \ref{conj2}
imply the Uniform and Quasi-uniform Resolution Conjectures (for the statement of
these conjectures see the Introduction).

\begin{prop}\label{conjqunif}
The SHGH Conjecture together with
Conjectures \ref{conj1} and \ref{conj2}
imply the Uniform and Quasi-uniform Resolution Conjectures.
\end{prop}
\noindent {\em Proof}. Since uniform implies quasi-uniform, let $Z$ be a
quasi-uniform point scheme, i.e. $Z=m\sum _{i=1,\dots ,9} P_i-\sum _{i=10,\dots
,n} m_iP_i$, $n\geq 9$, $m\geq m_1\geq \ldots m_n \geq 0$. We want to prove
that, assuming the SHGH Conjecture,
Conjecture \ref{conj1} and Conjecture \ref{conj2}, the map $\mu_k(Z)$, or
equivalently the map $\mu_F$ with $F=kL-m\sum _{i=1,\dots ,9} E_i-\sum
_{i=10,\dots ,n} m_iE_i$, is of maximal rank.
\par We can write $F=(k-3m)L-mK_X+\sum _{i\geq 10}(m-m_i) E_i$.
We can assume that $h_Z(k)>0$, otherwise $\mu_k(Z)$ is the zero map, hence
trivially injective;  since $Z$ is quasi-uniform, the SHGH conjecture then says
(see introduction) that $h_Z(k)= {k+2\choose 2}-9{m+1\choose2}
-\sum_i{m_i+1\choose2}$. In particular ${k+2\choose 2}-9{m+1\choose2}>0$, which
gives $k\geq 3m$. If $k=3m$, then $n=9$, in which case $F=m(3L-E_1-\cdots-E_9)$,
so $h^0(F)=1$ and $\mu_F$ has maximal rank.
\par Now let $k>3m$. In order to prove that $\mu_F$ has maximal rank, by
\ref{conj1} and \ref{conj2}
it is enough to prove that $\delta_0(C,Z,k)=0$ for each $\tilde C=dL-\sum
r_iE_i$ smooth rational in $X$; since
$\delta_0(C,Z,k)=\hbox{max}(0,\; -F\cdot \tilde C +a-1) + \hbox{max}(0,\;
-F\cdot \tilde C+b-1)$ with $a\leq b \leq d$, we'll just prove that $-F\cdot
\tilde C+b-1\leq 0$.
\par By the SHGH Conjecture, $\tilde C^2\geq -1$, so by adjunction formula
$K_X\cdot \tilde C=-\tilde C^2-2\leq-1$.   We hence find:$-F\cdot \tilde
C+b-1=(-(k-3m)L+mK_X-\sum _{i\geq 10}(m-m_i) E_i)\cdot \tilde C+b-1\leq
-(k-3m)d-m-\sum _{i\geq 10}r_i(m-m_i) +d-1<0$.

\prfend

\bigskip \par  {Alessandro Gimigliano, 
Dipartimento di Matematica and CIRAM,
Universit\`a di Bologna,
40126 Bologna, Italy  \quad email: gimiglia@dm.unibo.it

\medskip \par  Brian Harbourne,
Department of Mathematics,
University of Nebraska,
Lincoln, NE 68588-0323 USA  \quad email: bharbour@math.unl.edu

\medskip \par  Monica Id\`a,
Dipartimento di Matematica,
Universit\`a di Bologna,
40126 Bologna, Italy \quad email: ida@dm.unibo.it}


\begin{thebibliography} {biblio}

\bibitem[AH] {refAH} J. Alexander and A. Hirschowitz,  {\it
Polynomial interpolation
in several variables},
J.\ Alg.\ Geom. 4 (1995), 201-222.

\bibitem[As] {refAs} M.-G. Ascenzi, {\it The restricted tangent
bundle of a rational
curve in $\PP2$},  Comm. Algebra  16  (1988),  no. 11, 2193-2208.

\bibitem[BI] {refBI} E. Ballico and M. Id\`a, {\it On the minimal free
resolution for fat point schemes of multiplicity at most 3 in $\PP2$}, preprint
2006.

\bibitem[Ca] {refCa} M.V. Catalisano, {\it "Fat" points on a conic},
Comm. Algebra 19 (1991), 2153-2168.

\bibitem[CCMO] {refCCMO} C. Ciliberto, F. Cioffi, R. Miranda and F. Orecchia.
{\it Bivariate Hermite interpolation and linear systems of
plane curves with  base fat points}, in: Computer mathematics,
87-102, Lecture Notes Ser. Comput., 10, World Sci. Publishing,
River Edge, NJ,  2003.

\bibitem[E] {refEvain} L. Evain, {\it Computing limit linear series
with infinitesimal
methods}, preprint 2004 (arXiv:math.AG/0407143).

\bibitem[F1] {refF1} S. Fitchett, {\it  On Bounding the Number of
Generators for
Fat Point Ideals on the Projective Plane}, Journal of Algebra, 236
(2001), 502-521.

\bibitem[F2] {refF2} S. Fitchett, {\it Corrigendum to: "On bounding the
number of generators for fat point  ideals on the projective plane"
[J. Algebra 236 (2001), no. 2,  502--521]},  J. Algebra  276  (2004),
no. 1, 417-419.

\bibitem[F3] {refF3} S. Fitchett, {\it  Maps of linear systems on blow ups
of the Projective Plane}, J. Pure Appl. Algebra 156 (2001), 1-14.

\bibitem[FHH] {refFHH} S. Fitchett, B. Harbourne and S. Holay,
{\it Resolutions of Fat Point Ideals Involving Eight
General Points of $\PP2$}, J. Algebra 244 (2001), 684-705.

\bibitem[GM] {refGM} A.V. Geramita and P. Maroscia, {\it The ideal of
forms vanishing at a finite set of points in $\PPP ^n$}, J.Algebra 90
(1984), 528-555.

\bibitem[G] {refGi} A. Gimigliano, {\it On linear systems of
plane curves}, Thesis, Queen's University, Kingston (1987).

\bibitem[GHI] {refGHIpv} A. Gimigliano, B. Harbourne and M. Id\`a,
{\it Betti numbers for fat point ideals in the plane: a geometric approach},
preprint, http://www.math.unl.edu/~bharbour/shortBMS27-12-06Posted.pdf.

\bibitem[GI] {refGiId} A. Gimigliano and M. Id\`a, {\it The ideal resolution
for generic 3-fat points in $\PP2$},
J. Pure Appl. Algebra  187  (2004),  no. 1-3, 99-128.

\bibitem[GS] {refGS} D. Grayson, and M. Stillman, {\it Macaulay 2,
a software system for research in algebraic geometry},
Available at {\tt http://www.math.uiuc.edu/Macaulay2/}.

\bibitem[GH]{refGH} E. Guardo and B. Harbourne,
{\it Resolutions of ideals of six fat points in $\PP2$},
to appear, J. Alg. (arXiv:math.AG/0506611).

\bibitem[Ha1]{refSiena}  B. Harbourne, {\it The (unexpected)
importance of knowing $\alpha$}, pp. 267-272, in
``Projective Varieties with Unexpected Properties,
A Volume in Memory of Guiseppe Veronese,''
Proceedings of the international conference
{\it Varieties with Unexpected Properties}, Siena, Italy, June 8-13, 2004;
published 2005.


\bibitem[Ha2] {refIGP} B. Harbourne, {\it The Ideal Generation
Problem for Fat Points},
J. Pure Appl. Alg. 145(2), 165-182 (2000).

\bibitem[Ha3] {refanti} B. Harbourne, {\it Anticanonical rational surfaces},
Trans. Amer. Math. Soc. 349, 1191-1208 (1997).

\bibitem[Ha4] {refVanc} B. Harbourne, {\it The Geometry of rational
surfaces and Hilbert functions of points in the plane}. Can. Math.
Soc. Conf. Proc., vol. 6 (1986), 95-111.

\bibitem[Ha5]{refBHProc} B. Harbourne, {\it Rational Surfaces with $K^2>0$},
Proc. Amer. Math. Soc. 124, 727-733 (1996).

\bibitem[Ha6]{refCJM} B. Harbourne, {\it An Algorithm for Fat Points
on $\PP2$},
Can. J. Math. 52 (2000), 123-140.

\bibitem[HR] {refHRa} B. Harbourne and J. Ro\'e.
{\it Linear systems with multiple base points in ${\bf P}^2$},
Adv.\ Geom. 4 (2004), 41-59.

\bibitem[HHF] {refHHF} B. Harbourne, S. Holay and S. Fitchett,
{\it Resolutions of ideals of
quasiuniform fat point subschemes of ${\bf P}^2$},
Trans. Amer. Math. Soc. 355 (2003), no. 2, 593-608.

\bibitem[Ht] {refHt} R. Hartshorne, {\it Algebraic geometry}, Graduate Texts in
Mathematics, No. 52. Springer-Verlag, New York-Heidelberg,  1977.

\bibitem[Hi1] {refHi1}  A. Hirschowitz, {\it La M\'ethode d'Horace pour
l'interpolation \`a plusieurs variables}, Manuscripta Math. {\bf 50}
(1985), 337-388.

\bibitem[Hi2] {refHi} A. Hirschowitz,
{\it Une conjecture pour la cohomologie
des diviseurs sur les surfaces rationelles g\'en\'e\-ri\-ques},
Journ.\ Reine Angew.\ Math. 397
(1989), 208-213.

\bibitem[I1] {refId3} M. Id\`a, {\it On the homogeneous ideal of the
generic union
of lines in $\PP 3$}, J.Reine Angew.Math. 403 (1990), 67-153.

\bibitem[I2] {refId1} M. Id\`a,
{\it The minimal free resolution for the first infinitesimal
neighborhoods of $n$ general points in the plane},
J.\ Alg. 216 (1999), 741-753.

\bibitem[I3] {refId2} M. Id\`a,
{\it Generators for the generic rational space curve: low degree cases}.
Lecture Notes in Pure and Applied Math. Dekker 206 (1999), 169-210.

\bibitem[Mi] {refmig}
T. Mignon, {\it Syst\`emes de courbes planes \`a singularit\'es
impos\'ees: le cas des multiplicit\'es inf\'erieures ou \'egales
\`a quatre}, J.\ Pure Appl.\ Algebra 151 (2000), no. 2, 173-195.


\bibitem[Mu2] {refMu2} D. Mumford,  {\it Lectures on curves on an
algebraic surface},
Princeton 1966.

\bibitem[N] {refNtwo} M. Nagata, {\it On rational surfaces, II},
Mem.\ Coll.\ Sci.\
Univ.\ Kyoto, Ser.\ A Math.\ 33 (1960), 271-293.

\bibitem[OSS] {refOSS} C. Okonek, M. Schneider, H. Spindler, {\it Vector
bundles on
Complex Projective Spaces}, Progress in Mathematics, No. 3. Birkhauser, Boston,
Basel, Stuttgart 1980.

\bibitem[S] {refSe} B. Segre, {\it Alcune questioni su insiemi finiti
di punti in Geometria Algebrica},
Atti del Convegno Internaz. di Geom. Alg., Torino (1961).

\bibitem[Y] {refyang} S. Yang, {\it Linear systems in $\PP2$ with base points of
bounded multiplicity}, J. Algebraic Geometry 16 (2007), 19-38 (arXiv:math.AG/0406591).

\end{thebibliography}
\end{document}